# POWER OF THE SPACING TEST FOR LEAST-ANGLE REGRESSION

JEAN-MARC AZAÏS, YOHANN DE CASTRO, AND STÉPHANE MOURAREAU

ABSTRACT. Recent advances in Post-Selection Inference have shown that conditional testing is relevant and tractable in high-dimensions. In the Gaussian linear model, further works have derived unconditional test statistics such as the Kac-Rice Pivot for general penalized problems. In order to test the global null, a prominent offspring of this breakthrough is the spacing test that accounts the relative separation between the first two knots of the celebrated least-angle regression (LARS) algorithm. However, no results have been shown regarding the distribution of these test statistics under the alternative. For the first time, this paper addresses this important issue for the spacing test and shows that it is unconditionally unbiased. Furthermore, we provide the first extension of the spacing test to the frame of unknown noise variance.

More precisely, we investigate the power of the spacing test for LARS and prove that it is unbiased: its power is always greater or equal to the significance level $\alpha$. In particular, we describe the power of this test under various scenarii: we prove that its rejection region is optimal when the predictors are orthogonal; as the level $\alpha$ goes to zero, we show that the probability of getting a true positive is much greater than $\alpha$; and we give a detailed description of its power in the case of two predictors. Moreover, we numerically investigate a comparison between the spacing test for LARS and the Pearson's chi-squared test (goodness of fit).

When the noise variance is unknown, our analysis unleashes a new test statistic that can be computed in cubic time in the population size and which we refer to as the t-spacing test for LARS. The t-spacing test for LARS involves the first two knots of the LARS algorithm and we give its distribution under the null hypothesis. Interestingly, numerical experiments witness that the t-spacing test for LARS enjoys the same aforementioned properties as the spacing test for LARS.

## 1. INTRODUCTION

A major development in modern statistics has been brought by the idea that one can recover a high-dimensional target $\beta^\star$ from few linear observations $Y$ by $\ell_1$-minimization as soon as the target vector is "*sparse*" in a well-chosen basis. Undoubtedly, the notion of "*sparsity*" has encountered a large echo among the statistical community and many successful applications rely on $\ell_1$-minimization, the reader may consult [CDS98, Tib96, Fuc05, CT06, CT07] for some seminal works, [HTF09, BVDG11] for a review and references therein. More precisely, some of the most popular estimators in high-dimensional statistics remain the lasso [Tib96] and the Dantzig selector [CT07]. A large amount of interest has been dedicated to the estimation, prediction or support recovery problems using these estimators. This body of work has been developed around sufficient conditions on the design matrix $X$ (such that *Restricted Isometry Property* [CT06], *Restricted Eigenvalue* [BRT09], *Compatibility* [vdGB09, BVDG11], *Universal Distortion* [DC13, BLPR11], $\mathbf{H}_{s,1}$ [JN11], or *Irrepresentability* [Fuc05], to name but a few) that enclose the spectral properties of the design matrix on the set of (almost) sparse vectors.

Using one of these properties, one can exploits the Karush-Kuhn-Tucker conditions to get oracle inequalities or a control on the support recovery error.

Aside from those issues some recent works have been focused on hypothesis testing using penalized problems, see for instance [LSST13, LTTT14a, TLT14, TLTT14] and references therein. Compared to the sparse recovery problems, very little work has been done in statistical testing in high dimensions. As a matter of fact, one of the main difficulty is that there is no tractable distribution of sparse estimators (even under the aforementioned standard conditions of high-dimensional statistics). A successful approach is then to take into account the influence of each predictor in the regression problem. More precisely, some recent works in Post-Selection Inference have shown that the selection events can be explicitly expressed as closed convex polytopes depending simply on the signs and the indices of the nonzero coefficients of the solutions of standard procedures in high-dimensional statistics (typically the solutions of the lasso). Furthermore, an important advance has been brought by a useful parametrization of these convex polytopes under the Gaussian linear model, see for instance the book [HTW15]. In detection testing, this is done by the first two "*knots*" of the *least-angle regression algorithm* (LARS for short) which is intimately related to the dual program of the $\ell_1$-minimization problem, see [EHJT04] for example.

1.1. **Hypothesis testing using LARS.** The usual frame of the regression problems in high-dimensions is the following. Given an outcome vector $Y \in \mathbb{R}^n$, a matrix of predictor variables (or design matrix) $X \in \mathbb{R}^{n \times p}$ and a variance-covariance matrix $\Sigma$ such that
$$Y = X\beta^\star + \xi \quad \text{with} \quad \xi \sim \mathcal{N}_n(0, \Sigma),$$
we are concerned with testing whether $\beta^\star$ is equal to some known $\beta_0^\star$ or not. Notice that the response variable $Y$ does not depend directly on $\beta^\star$ but rather on $X\beta^\star$. We understand that a detection test may be interested in discerning between two hypothesis on the target vector, namely
$$\mathbb{H}_0 : \text{``}\beta^\star \in \beta_0^\star + \ker(X)\text{''} \quad \text{against} \quad \mathbb{H}_1 : \text{``}\beta^\star \notin \beta_0^\star + \ker(X)\text{''},$$
where $\ker(X)$ denotes the kernel of the design matrix $X$. It can be equivalently formulated (subtracting $X\beta_0^\star$) as a detection test whose null hypothesis is given by
$$\mathbb{H}_0 : \text{``}\beta^\star \in \ker(X)\text{''} \quad \text{against} \quad \mathbb{H}_1 : \text{``}\beta^\star \notin \ker(X)\text{''}.$$
To this end, we consider the vector of correlations
$$U := X^\top Y \sim \mathcal{N}_p(\mu^\star, R),$$
where $\mu^\star := X^\top X \beta^\star$ and $R := X^\top \Sigma X$. Observe that the hypotheses $\mathbb{H}_0$ and $\mathbb{H}_1$ can be equivalently written as

($\star$)
$$\mathbb{H}_0 : \text{``}\mu^\star = 0\text{''} \quad \text{against} \quad \mathbb{H}_1 : \text{``}\mu^\star \neq 0\text{''},$$

and remark that the knowledge of the noise variance-covariance matrix $\Sigma$ is equivalent to the knowledge of the correlations variance-covariance matrix $R$.

1.2. **The Spacing test for LARS.** The test statistic we are considering was introduced in a larger context of penalization problems by the pioneering works in [TLTT14, TLT14]. As mentioned by the authors of [TLT14], the general test statistic "may seem complicated". However, it can be greatly simplified in the frame of the standard regression problems under a very mild assumption, namely

(H)
$$\forall i \in [\![1, p]\!], \quad R_{ii} := X_i^\top \Sigma X_i = 1.$$

Note that this assumption is not very restrictive because the columns $X_i$ of $X$ can always be scaled to get (H). In this case, the entries of $\beta^\star$ are scaled but nor $\mathbb{H}_0$ neither $\mathbb{H}_1$ are changed. Hence, without loss of generality, we admit to invoke an innocuous normalization on the columns of the design matrix. Remark also that (H) is satisfied under the stronger assumption

(H Lasso) $$\Sigma = \mathrm{Id}_n \quad \text{and} \quad \forall i \in [\![1,p]\!], \quad \|X_i\|_2^2 = 1.$$

Moreover, observe that, almost surely, there exists a unique couple $(\hat\imath, \hat\varepsilon) \in [\![1,p]\!] \times \{\pm 1\}$ such that $\hat\varepsilon U_{\hat\imath} = \|U\|_\infty$. Under Assumption (H), the test statistic, refered to as *Spacing test for LARS*, simplifies to

(Pivot) $$S := \frac{\bar\Phi(\lambda_1)}{\bar\Phi(\lambda_2)},$$

with we denote by $\Phi$ the cumulative distribution function of the standard normal distribution, $\bar\Phi = 1 - \Phi$ its complement, $\lambda_1 := \hat\varepsilon U_{\hat\imath}$ the largest knot in the *lasso path* [EHJT04] and

$$\lambda_2 := \bigvee_{1 \leq j \neq \hat\imath \leq p} \left\{ \frac{U_j - R_{j\hat\imath} U_{\hat\imath}}{1 - \hat\varepsilon R_{j\hat\imath}} \vee \frac{-U_j + R_{j\hat\imath} U_{\hat\imath}}{1 + \hat\varepsilon R_{j\hat\imath}} \right\},$$

where $a \vee b := \max(a, b)$ and $U_i$ denotes the $i$-th entry of the vector $U$. Under Assumption (H Lasso), one has $R = X^\top X$ and $\lambda_2$ simplifies to the second largest knot in the *lasso path*. Interestingly, the authors of [TLT14] have shown that the test statistic $S$ is uniformly distributed on $[0,1]$ under the null hypothesis $\mathbb{H}_0$,

$$S \sim \mathrm{Unif}([0,1]).$$

Moreover, they derived the following rejection region

$$\mathrm{Reject}_\alpha := \{S \leq \alpha\},$$

for all $\alpha \in (0,1)$. In other words, the observed value of the test statistic $S$ is the $p$-value of the Spacing test for LARS.

Remark that the statistic $1-S$ is uniformly distributed on $[0,1]$, as well as many other transformations of the test statistic $S$. It may appear that the choice of rejection region $\mathrm{Reject}_\alpha$ is somehow arbitrary. Nevertheless, one can empirically witness (see Figure 1 for instance) that the *Spacing test for LARS* is an interesting test statistic that may take smaller values under the alternative hypothesis. However, no theoretical guarantees have been shown regarding its power. Furthermore, the *Spacing test for LARS* relies on the assumption that the variance-covariance matrix $\Sigma$ of the the noise is known and it should be interesting to bypass this limitation. To the best of our knowledge, this paper is the first to address these issues.

1.3. **Power of the Spacing test for LARS.** Recall that the Spacing test for LARS rejects $\mathbb{H}_0$ in favor of $\mathbb{H}_1$ when $\{S \leq \alpha\}$ occurs, where $S$ is defined by (Pivot). We assume that the noise variance-covariance matrix $\Sigma$ is known. We also assume that the columns $(X_i)_{i=1}^p$ of the design matrix $X$ are pairwise different and normalized with respect to Assumption (H). The first result shows that the Spacing test for LARS is unbiased.

**Theorem 1.** *Let $\alpha \in (0,1)$ be a significance level. Assume that the variance-covariance matrix $\Sigma$ of the noise is known and assume that Assumption (H) holds. Then, the Spacing test for LARS is unbiased: its power under the alternative is always greater or equal to the significance level $\alpha$.*

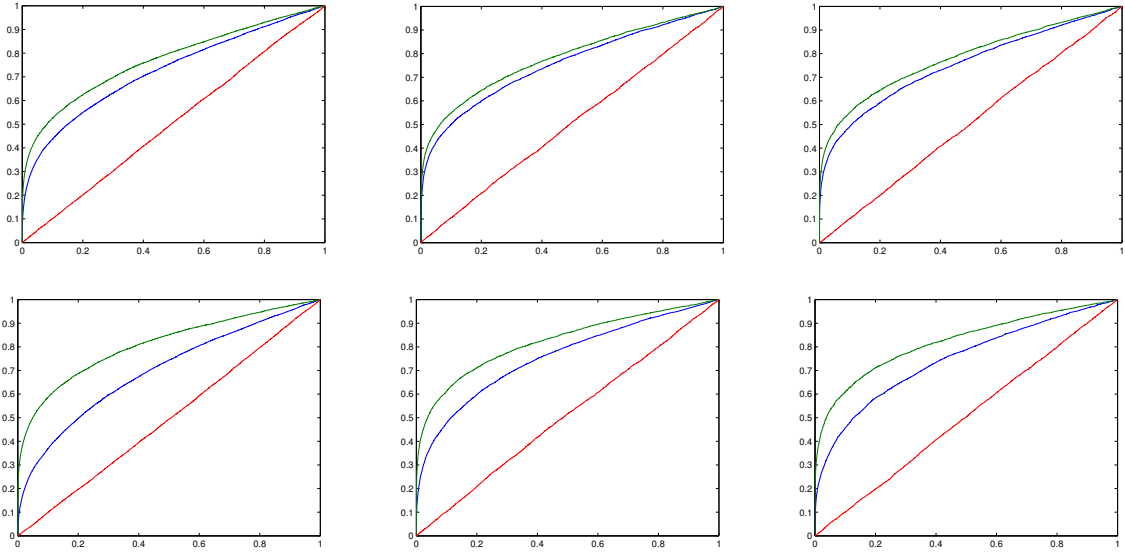

FIGURE 1. On each figure, empirical distribution function of 15,000 p-values coming from various scenarii. 5,000 p-values drawn under the null (red), 5,000 p-values of $S$ under the alternative (green) and 5,000 p-values of $T$ under the alternative (blue). At the top, the level of sparsity $s$ is equal to 2. At the bottom, $s$ is 5. In both cases, from left to right, $(n, p) = (50, 100)$, $(100, 200)$ and $(100, 500)$.

Under mild assumptions, this theorem ensures that the probability of getting a *true positive* is greater or equal to the probability of a *false positive*. Moreover, in the limit case when the significance level $\alpha$ goes to zero, this result is refined by Theorem 5: the probability of a *true positive* is much greater than the probability of getting a *false positive*. As a matter of fact, we prove that the cumulative distribution function of $S$ has a vertical tangent at the origin under the alternative hypothesis. The reader may consult Figure 1 which represents the empirical distribution function of $S$ that exactly describes the uniform law.

A proof of Theorem 1 can be found in Section 2.3. Interestingly, our proof is based on Anderson's inequality [And55] for symmetric convex sets. Moreover, we derive a simple and short proof of the distribution of the test statistic (Pivot) under the null, see Corollary 1 of Proposition 4.

Theorem 1 has a stronger version in the case of orthogonal designs, e.g. when the variance-covariance matrix $\Sigma$ is $\mathrm{Id}_n$ and $X^\top X = \mathrm{Id}_p$ (which implies that $n \geq p$).

**Theorem 2** (Orthogonal design). *Assume that $R = \mathrm{Id}_n$ then, under any alternative in $\mathbb{H}_1$, the density function of $S$ is decreasing. Hence, for all significance level $\alpha \in (0, 1)$, the region* $\mathrm{Reject}_\alpha = \{S \leq \alpha\}$ *is the most powerful region among all possible regions.*

This theorem may be seen as an evidence in favor of the choice of the rejection region as $\mathrm{Reject}_\alpha = \{S \leq \alpha\}$. A proof of Theorem 2 can be found in Section 2.4.

1.4. **Extension to unknown variance.** Interestingly, we can derive from our analysis a *studentization* of the test statistic (Pivot). Indeed, we consider the test statistic

(t-Pivot) $$T := \frac{1 - \mathbb{F}_{n-1}(T_1)}{1 - \mathbb{F}_{n-1}(T_2)},$$

where $\mathbb{F}_{n-1}$ denote the cumulative distribution function of the $t$-distribution with $n-1$ degrees of freedom and $T_1, T_2$ are statistics that can be computed in cubic time (cost of one Singular Value Decomposition (SVD) of the design matrix) from the first knots of the LARS algorithm, see Algorithm 1. In the sequel, for each $i \in [\![1,p]\!]$, we may denote by $X_{-i} \in \mathbb{R}^{n \times (p-1)}$ the sub-matrix of $X$ where the $i$-th column $X_i$ has been deleted and we may assume that it has rank $n$. Observe that this is a mild assumption in a high-dimensional context.

**Theorem 3** (t-Spacing test for LARS). *Assume that the variance-covariance matrix $\Sigma$ is $\sigma^2 \mathrm{Id}_n$ where $\sigma > 0$ is unknown and that for all $i \neq j \in [\![1,p]\!]$, one has $\|X_i\|_2 = 1$, $X_i \neq X_j$ and $X_{-i}$ has rank $n$. Then, under the null $\mathbb{H}_0$, the statistic $T$ described by Algorithm 1 is uniformly distributed on $[0,1]$.*

In particular, we derive a detection test of significance level $\alpha$ considering the rejection region $\mathrm{Reject}_\alpha = \{T \leq \alpha\}$. A proof of Theorem 3 can be found in Section 3. One can empirically witness (see Figure 1 for instance) that the *t-Spacing test for LARS* is an interesting test statistic that may take smaller values under the alternative hypothesis.

---

**Algorithm 1:** t-Spacing test

**Data**: An observation $Y \in \mathbb{R}^n$ and a design matrix $X \in \mathbb{R}^{n \times p}$.
**Result**: A $p$-value $T \in (0,1)$.

*Compute the first LARS knot $\lambda_1$;*

(1) Set $U := X^\top Y$;
(2) Find $(\hat{\imath}, \hat{\varepsilon}) \in [\![1,p]\!] \times \{\pm 1\}$ such that $\hat{\varepsilon} U_{\hat{\imath}} = \|U\|_\infty$ and set $\lambda_1 := \hat{\varepsilon} U_{\hat{\imath}}$;

*Compute the second LARS knot $\lambda_2$;*

(3) Set $R := X^\top X$;
(4) Set
$$\lambda_2 := \bigvee_{1 \leq j \neq \hat{\imath} \leq p} \left\{ \frac{U_j - R_{j\hat{\imath}} U_{\hat{\imath}}}{1 - \hat{\varepsilon} R_{j\hat{\imath}}} \vee \frac{-U_j + R_{j\hat{\imath}} U_{\hat{\imath}}}{1 + \hat{\varepsilon} R_{j\hat{\imath}}} \right\};$$

*Compute the variance estimator $\hat{\sigma}$;*

(5) Set $R_{-\hat{\imath}} := X_{-\hat{\imath}}^\top (\mathrm{Id}_n - X_{\hat{\imath}} X_{\hat{\imath}}^\top) X_{-\hat{\imath}}$;
(6) Compute $R_{-\hat{\imath}}^{-1/2}$ the square root of the pseudoinverse of $R_{-\hat{\imath}}$;
(7) Set
$$\hat{\sigma} := \frac{\|R_{-\hat{\imath}}^{-1/2} V_{-\hat{\imath}}\|_2}{\sqrt{n-1}},$$
where
$$V_{-\hat{\imath}} := (U_1 - R_{1\hat{\imath}} U_{\hat{\imath}}, \ldots, U_{\hat{\imath}-1} - R_{(\hat{\imath}-1)\hat{\imath}} U_{\hat{\imath}}, U_{\hat{\imath}+1} - R_{(\hat{\imath}+1)\hat{\imath}} U_{\hat{\imath}}, \ldots, U_p - R_{p\hat{\imath}} U_{\hat{\imath}});$$

*Compute the p-value $T$;*

(8) Set $T_1 := \lambda_1/\hat{\sigma}$ and $T_2 := \lambda_2/\hat{\sigma}$;
(9) Set
$$T := \frac{1 - \mathbb{F}_{n-1}(T_1)}{1 - \mathbb{F}_{n-1}(T_2)},$$
where we denote by $\mathbb{F}_{n-1}$ the cumulative distribution function of the $t$-distribution with $n-1$ degree(s) of freedom.

Observe that Algorithm 1 requires the computation of one SVD at step 6. We deduce that its computational cost is $\mathcal{O}(p^3)$ which is reasonable in high-dimensional statistics.

1.5. **Empirical distributions of the p-values.** Figure 1 shows the empirical distribution of a sample of 15,000 p-values constructed from standard regression problems under the global null and under the alternative for the pivots $S$ and $T$. Design matrices $X$ and the mean (under the alternative) have been drawn uniformly at random from the following cases

- $X$ is a design matrix of size $50 \times 100$, $100 \times 200$ or $100 \times 500$ with i.i.d. $\mathcal{N}(0,1)$ entries.
- $\beta$ is a vector with i.i.d. $\mathcal{N}(0,1)$ (small mean), $\mathcal{N}(0,4)$ (medium mean) or $\mathcal{N}(\sqrt{2\log p}, 1)$ (high mean) entries.

Under the null, the agreement with uniform is very strong. Moreover, the Spacing test for LARS is empirically more powerful than the *t-Spacing test for LARS* and both seem to be unbiased. However, in a context of very high-dimensional regression, the *t-Spacing test for LARS* is very similar to the *Spacing test for LARS* due to standard results on Student and chi-squared distribution.

1.6. **Previous works.** Our test can be also referred to as the *Kac-Rice test* as introduced in the broader frame of penalization problems in the seminal paper [TLT14]. The interested reader may consult Theorem 1 in [TLT14] where the general "*Kac-Rice pivot*" is defined. Note that various important results on this subject have been obtained recently and we do not pursue on a comprehensive study here. The interested reader may consult Chapter 6 of the captivating book [HTW15].

The statistic *Kac-Rice pivot* given in [TLT14] has been used for model selection and confidence intervals on the target entries. In the frame of lasso, the optimality of these approaches is discussed in [LSST13, TLTT14]. Interestingly, the Spacing test is a nonasymptotic version of the covariance test [LTTT14a, TLTT14], and is asymptotically equivalent to it. Note they have been intensively commented among the literature, see [LTTT14c, LTTT14b, BMvdG14] for instance.

1.7. **Organization of the paper.** The next section is devoted to the proof of the main results on the power. In particular, the reader may find the exact formulation of Theorem 5 mentioned in the introduction. Section 3 addresses the issue of extending the Spacing test for LARS to the unknown variance frame. Section 4 presents a fine description of the Spacing test for LARS' power in the case of two predictors. The last section gives a numerical comparison with the Pearson's chi-squared test (goodness of fit).

## 2. Power of the Spacing test for LARS

2.1. **Model and notation.** Recall that the vector of correlations $U = X^\top Y$ enjoys

$$U = (U_1, \ldots, U_p) \sim \begin{cases} \mathcal{N}_p(0, R) & \text{under the null hypothesis,} \\ \mathcal{N}_p(\mu^\star, R) & \text{under the alternative hypothesis,} \end{cases}$$

where $R = X^\top \Sigma X$ and $\mu^\star = X^\top X \beta^\star$. Indeed, observe that

$$\{\mu^\star = 0\} \Leftrightarrow \{\mathbb{H}_0 : \text{``}\beta^\star \in \ker(X)\text{''}\}.$$

It is well known (see for instance the book [HTW15]) that the first knot $\lambda_1$ of the LARS algorithm enjoys $\lambda_1 = \|U\|_\infty$. Assume that the columns of $X$ are pairwise different.

It implies that, with probability one, there exists a unique pair $(\hat{\imath}, \hat{\varepsilon})$ with $\hat{\imath} \in [\![1, p]\!]$, $\hat{\varepsilon} = \pm 1$ and such that

(1) $$\hat{\varepsilon} U_{\hat{\imath}} = \|U\|_\infty.$$

Observe that the events $\mathscr{E}_{i,\varepsilon} := \{\varepsilon U_i = \|U\|_\infty\}$ are almost surely disjoints, and note that

$$\lambda_1 = \sum_{i=1}^{p} \sum_{\varepsilon = \pm 1} \varepsilon U_i \mathbb{1}_{\mathscr{E}_{i,\varepsilon}},$$

where $\mathbb{1}$ denotes the indicator function. Write, for all $(i,j) \in [\![1,p]\!]^2$, $U_j = R_{ji} U_i + U_j^i$, the regression of $U_j$ onto $U_i$. Recall that the residuals $U_j^i$ are independent of $U_i$. Denote, for all $i \in [\![1, p]\!]$ and $\varepsilon = \pm 1$,

$$\lambda_2^{i,\varepsilon} = \bigvee_{1 \leq j \neq i \leq p} \left\{ \frac{U_j^i}{1 - \varepsilon R_{ji}} \vee \frac{-U_j^i}{1 + \varepsilon R_{ji}} \right\}.$$

Furthermore, remark that $\mathscr{E}_{i,\varepsilon} = \{\lambda_2^{i,\varepsilon} < \varepsilon U_i\}$. Indeed, for all $i \neq j \in [\![1, p]\!]$,

$$\{-\varepsilon U_i < U_j < \varepsilon U_i\} = \{-\varepsilon U_i (1 + \varepsilon R_{ji}) < U_j - R_{ji} U_i < \varepsilon U_i (1 - \varepsilon R_{ji})\},$$

$$= \left\{ \left\{ \frac{U_j^i}{1 - \varepsilon R_{ji}} \vee \frac{-U_j^i}{1 + \varepsilon R_{ji}} \right\} < \varepsilon U_i \right\}.$$

Hence, define the random variable $\lambda_2$ as

$$\lambda_2 = \sum_{i=1}^{p} \sum_{\varepsilon = \pm 1} \lambda_2^{i,\varepsilon} \mathbb{1}_{\mathscr{E}_{i,\varepsilon}}.$$

We deduce that

(2) $$(\lambda_1, \lambda_2) = \sum_{i=1}^{p} \sum_{\varepsilon = \pm 1} (\varepsilon U_i, \lambda_2^{i,\varepsilon}) \mathbb{1}_{\{\varepsilon U_i > \lambda_2^{i,\varepsilon}\}}.$$

Denote by $\varphi$ the probability density function of the standard normal distribution.

**Lemma 1.** *For each $i \in [\![1, p]\!]$ and $\varepsilon = \pm 1$, the random variable $\lambda_2^{i,\varepsilon}$ has a density $p_{\lambda_2^{i,\varepsilon}}^{\mu^\star}$. The joint density of $(\lambda_1, \lambda_2)$ is given by*

(3) $$\forall (\ell_1, \ell_2) \in \mathbb{R}^2, \quad p_{(\lambda_1, \lambda_2)}^{\mu^\star}(\ell_1, \ell_2) = \sum_{i=1}^{p} \sum_{\varepsilon = \pm 1} \varphi(\ell_1 - \varepsilon \mu_i^\star) p_{\lambda_2^{i,\varepsilon}}^{\mu^\star}(\ell_2) \mathbb{1}_{\{0 \leq \ell_2 \leq \ell_1\}}.$$

*Proof.* One can check that $\lambda_2^{i,\varepsilon}$ has a density, the reader may also consult Ylvisaker's theorem, see Theorem 1.22 in [AW09] for example.

Observe that for all $(i, \varepsilon) \in [\![1, p]\!] \times \{\pm 1\}$, the random variable $\lambda_2^{i,\varepsilon}$ is a deterministic function of the random variables $U_j^i$ for $j \neq i$ and hence it is independent of $\varepsilon U_i$. We get that the density function $p_{(\varepsilon U_i, \lambda_2^{i,\varepsilon})}^{\mu^\star}$ of $(\varepsilon U_i, \lambda_2^{i,\varepsilon})$ with respect to Lebesgue measure is given by

$$\forall (i, \varepsilon) \in [\![1, p]\!] \times \{\pm 1\}, \quad \forall (\ell_1, \ell_2) \in \mathbb{R}^2, \quad p_{(\varepsilon U_i, \lambda_2^{i,\varepsilon})}^{\mu^\star}(\ell_1, \ell_2) = \varphi(\ell_1 - \varepsilon \mu_i^\star) p_{\lambda_2^{i,\varepsilon}}^{\mu^\star}(\ell_2).$$

Invoke (2) to complete the proof. □

**Lemma 2.** *For the study the distribution of $S$, we can assume, without loss of generality, that the expectations $\mu_i^\star$ are non-negative.*

*Proof.* Let $\mu^\star \in \mathbb{R}^p$ and consider the linear map $T : \mathbb{R}^p \to \mathbb{R}^p$ that changes the signs of the coordinates of $U$ with negative expectation. Set

$$\bar{U} := T(U) = \{t_i U_i : i \in [\![1, p]\!]\},$$

where for all $i \in [\![1, p]\!]$, $t_i$ is the sign of $\mu^\star_i$. Each coordinate of $\bar{U}$ has non-negative expectation and the variance-covariance matrix of $\bar{U}$ is now $\bar{R}$ with $\bar{R}_{i,j} = t_i t_j R_{i,j}$.

Let us check, with obvious notation, that the test statistic $S$ enjoys $S(\bar{U}) = S(U)$. Indeed, it holds that the first knot $\lambda_1$ satisfies $\lambda_1(\bar{U}) = \lambda_1(U)$, $\mathcal{E}_{i,\varepsilon}(\bar{U}) = \mathcal{E}_{i,t_i\varepsilon}(U)$ and one can note that

$$\bar{U}^i_j = t_j U_j - \bar{R}_{ij} t_i U_i = t_j U_j - t_i t_j R_{ij} t_i U_i = t_j U^i_j,$$

and

$$\lambda_2^{i,\varepsilon}(\bar{U}) = \bigvee_{j \neq i} W_{i,\varepsilon,j}(\bar{U}) \quad \text{with} \quad W_{i,\varepsilon,j}(\bar{U}) := \frac{\varepsilon t_j U^i_j}{1 - t_i t_j R_{ij}} \vee \frac{-\varepsilon t_j U^i_j}{1 + t_i t_j R_{ij}}.$$

One may check that, whatever the signs $t_i, t_j$ are, it holds $W_{i,\varepsilon,j}(\bar{U}) = W_{i,(t_i\varepsilon),j}(U)$. Thus $\lambda_2^{i,\varepsilon}(\bar{U}) = \lambda_2^{i,(t_i\varepsilon)}(U)$ implying $\lambda_2(\bar{U}) = \lambda_2(U)$. □

2.2. **Piecewise calculus of the power.** We have the following useful proposition giving an exact expression of the power of Spacing test for LARS as weigthed sum of Gaussian mesures of disjoint cones. Denote by $\mathcal{C}_{i,\varepsilon}$ the cone

$$\mathcal{C}_{i,\varepsilon} := \{(u_1, \ldots, u_p) \in \mathbb{R}^p \quad \text{such that} \quad \forall j \neq i, \ |u_j| < \varepsilon u_i\},$$

recall that $\bar{\Phi} = 1 - \Phi$ is the complement of the standard normal cumulative distribution function and define by $\bar{\Phi}^{-1}$ its inverse function.

**Proposition 4.** *For all $\alpha \in (0, 1)$, define*

(4) $$h_\alpha(\ell) := \bar{\Phi}^{-1}(\alpha \bar{\Phi}(\ell)) - \ell,$$

*Then it holds,*

(5) $$\mathbb{P}_{\mu^\star}\{S \leq \alpha\} = \alpha \mathbb{E}_{\mu^\star}\Big\{\sum_{i=1}^p \sum_{\varepsilon=\pm 1} \exp[\varepsilon \mu^\star_i h_\alpha(\varepsilon U_i)] \mathbb{1}_{\{U \in \mathcal{C}_{i,\varepsilon}\}}\Big\},$$

*where $\mathbb{E}_{\mu^\star}$ denotes the expectation under the Gaussian distribution $\mathcal{N}_p(\mu^\star, R)$.*

*Proof.* Let $\alpha \in (0, 1)$. Note that

(6) $$\{S \leq \alpha\} = \{\lambda_1 \geq \bar{\Phi}^{-1}(\alpha/2)\} \cap \{\lambda_2 \leq \bar{\Phi}^{-1}(\bar{\Phi}(\lambda_1)/\alpha)\}.$$

Using (2), the change of variable $q_1 = \bar\Phi^{-1}(\bar\Phi(\ell_1)/\alpha)$ and (3), it holds

$$\mathbb{P}_{\mu^\star}\{S \leq \alpha\} = \sum_{i=1}^p \sum_{\varepsilon=\pm 1} \int_{\bar\Phi^{-1}(\frac{\alpha}{2})}^{+\infty} d\ell_1 \varphi(\ell_1 - \varepsilon\mu_i^\star) \int_0^{\bar\Phi^{-1}(\bar\Phi(\ell_1))/\alpha} d\ell_2\, p_{\lambda_2^{i,\varepsilon}}^{\mu^\star}(\ell_2),$$

$$= \alpha \sum_{i=1}^p \sum_{\varepsilon=\pm 1} \int_0^{+\infty} dq_1 \frac{\varphi(q_1)}{\varphi(\ell_1)} \varphi(\ell_1 - \varepsilon\mu_i^\star) \int_0^{q_1} d\ell_2\, p_{\lambda_2^{i,\varepsilon}}^{\mu^\star}(\ell_2),$$

$$= \alpha \sum_{i=1}^p \sum_{\varepsilon=\pm 1} \int_0^{+\infty} dq_1 e^{\varepsilon\mu_i^\star(\ell_1 - q_1)} \varphi(q_1 - \varepsilon\mu_i^\star) \int_0^{q_1} d\ell_2\, p_{\lambda_2^{i,\varepsilon}}^{\mu^\star}(\ell_2),$$

$$= \alpha \sum_{i=1}^p \sum_{\varepsilon=\pm 1} \int_0^{+\infty} dq_1 e^{\varepsilon\mu_i^\star(\ell_1 - q_1)} \int_0^{q_1} d\ell_2\, p_{(\varepsilon U_i, \lambda_2^{i,\varepsilon})}^{\mu^\star}(q_1, \ell_2),$$

$$= \alpha \sum_{i=1}^p \sum_{\varepsilon=\pm 1} \mathbb{E}_{\mu^\star}\left[ \exp\left[\varepsilon\mu_i^\star(\bar\Phi^{-1}(\alpha\bar\Phi(\varepsilon U_i))-\varepsilon U_i)\right] \mathbb{1}_{\{\varepsilon U_i > \lambda_2^{i,\varepsilon}\}} \right],$$

$$= \alpha \sum_{i=1}^p \sum_{\varepsilon=\pm 1} \mathbb{E}_{\mu^\star}\left[ \exp\left[\varepsilon\mu_i^\star(\bar\Phi^{-1}(\alpha\bar\Phi(\varepsilon U_i))-\varepsilon U_i)\right] \mathbb{1}_{\{U \in \mathscr{C}_{i,\varepsilon}\}} \right].$$

as claimed. □

*Remark.* Note the numerical evaluation of (5) can be performed using a $n$-dimensional integral, see Section 5.

**Corollary 1.** *Under $\mathbb{H}_0$, the statistics S defined by* (Pivot) *follows a uniform distribution on* $[0,1]$.

*Proof.* The null hypothesis is equivalent to $\mu^\star = 0$ and, from (5), we recover that

$$\mathbb{P}_0\{S \leq \alpha\} = \alpha \sum_{i=1}^p \sum_{\varepsilon=\pm 1} \mathbb{E}_0(\mathscr{C}_{i,\varepsilon}) = \alpha,$$

i.e. the level of Spacing test for LARS is $\alpha$. This proves that, under $\mathbb{H}_0$, the test statistics $S$ satisfies

$$S \sim \mathrm{Unif}([0,1]),$$

as claimed. □

2.3. **Distribution under the alternative.** This section is devoted to the proof of Theorem 1.

**Step 1:** By a standard approximation argument, one may assume that $R$ is a regular matrix. Indeed, if $R$ is singular we can approximate it by a sequence $(R_m)_{m \geq 0}$ of regular matrices with bounded variance. If for each of these matrices we have $\mathbb{P}_{\mu^\star}\{S \leq \alpha\} \geq \alpha$ then the result will pass to $R$ by dominated convergence in (5). Furthermore, using Lemma 2, we may also assume that $\forall i \in [\![1,p]\!]$, $\mu_i^\star \geq 0$.

Recall that $\mathscr{C}_{i,\varepsilon}$ is the cone

$$\mathscr{C}_{i,\varepsilon} := \{(u_1, \ldots, u_p) \in \mathbb{R}^p \quad \text{such that} \quad \forall j \neq i,\ |u_j| < \varepsilon u_i\}.$$

and denote by $\gamma$ the non-degenerate Gaussian measure associated with the multivariate normal distribution $\mathcal{N}_p(0, R)$.

**Step 2:** We start from (5) to get that

$$\frac{1}{\alpha}\mathbb{P}_{\mu^\star}\{S \leq \alpha\} = \mathbb{E}_{\mu^\star}\Big\{\sum_{i=1}^{p}\sum_{\varepsilon=\pm 1}\exp\big[\varepsilon\mu_i^\star h_\alpha(\varepsilon U_i)\big]\mathbb{1}_{\mathscr{C}_{i,\varepsilon}}\Big\},$$

$$\geq 1 + \mathbb{E}_{\mu^\star}\Big\{\sum_{i=1}^{p}\sum_{\varepsilon=\pm 1}\big[\varepsilon\mu_i^\star h_\alpha(\varepsilon U_i)\big]\mathbb{1}_{\mathscr{C}_{i,\varepsilon}}\Big\}.$$

Perform an integration using the fibers $F_{\ell,i,\varepsilon} := \{u_i = \varepsilon\ell\}\bigcap \mathscr{C}_{i,\varepsilon}$ to obtain that

$$\frac{1}{\alpha}\mathbb{P}_{\mu^\star}\{S \leq \alpha\} \geq 1 + \int_0^{+\infty}\sum_{i=1}^{p}\sum_{\varepsilon=\pm 1}\varepsilon\mu_i^\star h_\alpha(\ell)\sigma_{\mu^\star}(\ell,i,\varepsilon)d\ell,$$

where $\sigma_{\mu^\star}(\ell,i,\varepsilon)$ is the integral of the density function $\varphi_{\mu^\star}$ of the multivariate normal distribution $\mathscr{N}_p(\mu^\star,R)$ on the fiber $F_{\ell,i,\varepsilon}$.

**Step 3:** Let $\ell > 0$ and $a \geq 0$. Consider the hypercube $H_\ell := [-\ell,\ell]^p$ and denote by $H_\ell - a\mu^\star$ its translation by vector $-a\mu^\star$. Invoke Anderson's inequality (see Lemma 4) to get that

$$a \mapsto \gamma(H_\ell - a\mu^\star) := \mathbb{P}\{\mathscr{N}_p(0,R) \in H_\ell - a\mu^\star\},$$

is a non-increasing function on the domain $a \geq 0$. In particular, its derivative at point $a = 1$ is non-positive. It reads as

$$\lim_{\eta \to 0}\frac{1}{\eta}\big(\gamma(H_\ell - (1+\eta)\mu^\star) - \gamma(H_\ell - \mu^\star)\big) \leq 0,$$

and this quantity is simply, by **Step 4**,

(7) $$\sum_{i=1}^{p}\mu_i^\star\sigma_{\mu^\star}(\ell,i,-1) - \sum_{i=1}^{p}\mu_i^\star\sigma_{\mu^\star}(\ell,i,+1) \leq 0.$$

Finally, the positivity of $h_\alpha(\ell)$ (see Lemma 5) completes the proof.

**Step 4:** In the context of **Step 3**, computation on $\gamma(H_\ell - a\mu^\star)$ gives that

$$\frac{d}{da}\gamma(H_\ell - a\mu^\star) = \int_{H_\ell}\frac{\partial}{\partial a}\varphi_{a\mu^\star}(z)dz = \sum_{i=1}^{p}-\mu_i\int_{H_\ell}\frac{\partial}{\partial z_i}\varphi_{a\mu^\star}(z)dz$$

$$= \sum_{i=1}^{p}\sum_{\varepsilon=\pm 1}-\varepsilon\mu_i^\star\sigma_{a\mu^\star}(\ell,i,\varepsilon).$$

where, for all $a > 0$, we denote by $\sigma_{a\mu^\star}(\ell,i,\varepsilon)$ the integral of the density function $\varphi_{a\mu^\star}$ of the multivariate normal distribution $\mathscr{N}_p(a\mu^\star,R)$ on the fiber $F_{\ell,i,\varepsilon}$.

This computation might also be illustrated via finite difference method, one may see Figure 2 for instance.

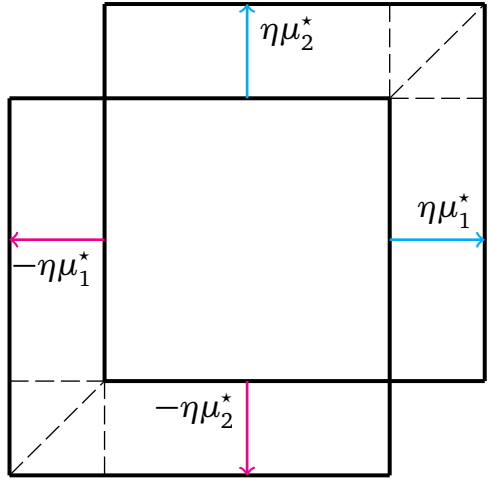

FIGURE 2. Illustration of (7) in dimension 2. Passing to the limit, contribution of triangles (dashed lines) vanish and the derivative in $a = 1$ is equal to the sum of each face with a weight $\varepsilon \mu_i^\star$ corresponding to its orientation.

2.4. **Orthogonal case.** In this section, we give the proof of Theorem 2. Invoke (3) to get that, under $\mathbb{H}_1$,

$$\forall (\ell_1, \ell_2) \in \mathbb{R}^2, \quad p_{(\lambda_1, \lambda_2)}^{\mu^\star}(\ell_1, \ell_2) = \sum_{i=1}^{p} \sum_{\varepsilon = \pm 1} \varphi(\ell_1 - \varepsilon \mu_i^\star) p_{\lambda_2^{i,\varepsilon}}^{\mu^\star}(\ell_2) \mathbb{1}_{\{0 \leq \ell_2 \leq \ell_1\}}.$$

Recall that $\hat{\imath} \in [\![1, p]\!]$ is defined by (1). Since $R = \mathrm{Id}_p$, remark that

$$\lambda_2 = \max_{j \neq \hat{\imath}} |U_j|,$$

Furthermore, observe that $\lambda_2^{i,+1} = \lambda_2^{i,-1}$ almost surely. It implies that for all $i \in [\![1, p]\!]$,

$$p_{\lambda_2^{i,+1}}^{\mu^\star} = p_{\lambda_2^{i,-1}}^{\mu^\star}.$$

Denote by $p_{\lambda_2^i}^{\mu^\star}$ their common value. As a consequence, it holds

$$\forall (\ell_1, \ell_2) \in \mathbb{R}^2, \quad p_{(\lambda_1, \lambda_2)}^{\mu^\star}(\ell_1, \ell_2) = \sum_{i=1}^{p} \big(\varphi(\ell_1 - \mu_i) + \varphi(\ell_1 + \mu_i)\big) p_{\lambda_2^i}^{\mu^\star}(\ell_2) \mathbb{1}_{\{0 \leq \ell_2 \leq \ell_1\}}.$$

It implies that, conditionally to $\lambda_2 = \ell_2$, the random variable $\bar{\Phi}(\lambda_1)$ admits the density

$$(8) \quad p_{(\bar{\Phi}(\lambda_1) | \lambda_2 = \ell_2)}(v) = (const) \sum_{i=1}^{p} \cosh\big(\bar{\Phi}^{-1}(v) \mu_i^\star\big) \mathbb{1}_{\{\bar{\Phi}^{-1}(v) \geq \ell_2\}}.$$

Since $\bar{\Phi}^{-1}(v)$ remains in the positive domain, the functions into the sum above are non-increasing and strictly decreasing for the index $i$ such that $\mu_i^\star > 0$. We have clearly the same result for the expression equivalent to (8) given the conditional density of $S$.

Deconditionning we obtain that the density of $S$ is a mixture of non-increasing functions, thus non-increasing. In addition the deconditioning formula gives positive weights to decreasing functions thus, in fact, the density is decreasing.

## 2.5. Asymptotic case.

**Theorem 5.** *Under* $\mathbb{H}_1$, *it holds*

$$\alpha^{-1}\mathbb{P}_{\mu^\star}\{S \leq \alpha\} \to +\infty,$$

*as $\alpha$ goes to zero, where $\mu^\star = (X^\top X)\beta^\star$ and $\mathbb{P}_{\mu^\star}$ denotes the law of $X^\top Y \sim \mathcal{N}_p(\mu^\star, R)$.*

*Proof.* Recall that $\mathbb{H}_1$ is equivalent to $\mu^\star \neq 0$. Without loss of generality, assume $\mu_1^\star > 0$ and note that

$$\alpha^{-1}\mathbb{P}_{\mu^\star}\{S \leq \alpha\} = \mathbb{E}_{\mu^\star}\Big\{\sum_{i=1}^p \sum_{\varepsilon=\pm 1} \exp\big[\mu_i^\star h_\alpha(\varepsilon U_i)\big]\mathbb{1}_{\{U \in \mathscr{C}_{i,\varepsilon}\}}\Big\}$$

$$\geq \mathbb{E}_{\mu^\star}\Big\{\exp\big[\mu_1^\star h_\alpha(U_1)\big]\mathbb{1}_{\{U \in \mathscr{C}_{1,1}\}}\Big\} =: \Lambda_{1,1}(\alpha).$$

Moreover, observe that

$$\forall x \in \mathbb{R}, \quad h_\alpha(x) \to +\infty,$$

as $\alpha > 0$ goes to zero. In particular, it yields

$$\forall x \in \mathbb{R}, \quad \exp\big[\mu_1^\star h_\alpha(x)\big] \to +\infty,$$

as $\alpha > 0$ goes to zero. Eventually, let $(\alpha_n)_{n\in\mathbb{N}}$ be any sequence of positive reals that goes to zero as $n$ tends to $\infty$. Invoke Fatou's lemma to get that

$$\lim_{n\to\infty} \alpha_n^{-1}\mathbb{P}_{\mu^\star}\{S \leq \alpha_n\} \geq \lim_{n\to\infty} \Lambda_{1,1}(\alpha_n)$$

$$\geq \liminf_{n\to\infty} \mathbb{E}_{\mu^\star}\Big\{\exp\big[\mu_1^\star h_\alpha(U_1)\big]\mathbb{1}_{\{U \in \mathscr{C}_{1,1}\}}\Big\}$$

$$\geq \mathbb{E}_{\mu^\star}\Big\{\liminf_{\alpha_n\to 0}\exp\big[\mu_1^\star h_\alpha(U_1)\big]\mathbb{1}_{\{U \in \mathscr{C}_{1,1}\}}\Big\}$$

which concludes the proof. □

## 3. Studentization of the Spacing test for LARS

In this section, we give the proof of Theorem 3.

### 3.1. Model and notation.
Assume that the variance-covariance matrix $\Sigma$ of the Gaussian noise $\xi$ is $\sigma^2\mathrm{Id}_n$ where $\sigma > 0$ is unknown. Assume also that the columns $(X_i)_{i=1}^p$ of the design matrix $X$ enjoy $\|X_i\|_2 = 1$ and denote by $U := X^\top Y$ the correlation vector satisfying

$$U = (U_1, \ldots, U_p) \sim \begin{cases} \mathcal{N}_p(0, \sigma^2 R) & \text{under the null hypothesis,} \\ \mathcal{N}_p(\mu^\star, \sigma^2 R) & \text{under the alternative hypothesis,} \end{cases}$$

where $R = X^\top X$ and $\mu^\star = R\beta^\star$. Observe that the knots of the LARS algorithm are given by

$$(\lambda_1, \lambda_2) = \sum_{i=1}^p \sum_{\varepsilon=\pm 1} (\varepsilon U_i, \lambda_2^{i,\varepsilon})\mathbb{1}_{\{\varepsilon U_i > \lambda_2^{i,\varepsilon}\}}.$$

For each $i \in [\![1,p]\!]$, we denote by $X_{-i} \in \mathbb{R}^{n\times(p-1)}$ the sub-matrix of $X$ where the $i$-th column $X_i$ has been deleted. Also, we denote by $U_{-i} \in \mathbb{R}^{p-1}$ (resp. $\mu_{-i}^\star$) the sub-vector of $U$ (resp. $\mu^\star$) where the $i$-th entry has been deleted. Observe that the regression of $U_{-i}$ onto $U_i$ reads

$$U_{-i} = (R_i)_{-i}U_i + V_{-i},$$

where $(R_i)_{-i} \in \mathbb{R}^{p-1}$ denotes the sub-vector of the $i$-th column $R_i$ of the matrix $R$ where the $i$-th entry has been deleted. Observe that the vector $V_{-i} \in \mathbb{R}^{p-1}$ is a Gaussian vector independent of $U_i$ such that

$$(9) \qquad V_{-i} \sim \mathcal{N}_{p-1}(\mu_{-i}^\star - (R_i)_{-i}\mu_i^\star, \sigma^2 R_{-i}),$$

where $R_{-i} := X_{-i}^\top(\mathrm{Id}_n - X_i X_i^\top)X_{-i}$ denotes its variance-covariance matrix. Notice that if $X_{-i}$ has full rank (namely $n$) then $R_{-i}$ has rank $n-1$. Denote $R_{-i}^{-1/2}$ the only symmetric matrix such that $R_{-i}^{-1/2} R_{-i} R_{-i}^{-1/2}$ is the orthogonal projection onto the range of $R_{-i}$ (observe that $R_{-i}^{-1/2}$ is the square root of the Moore-Penrose pseudoinverse of $R_{-i}$).

## 3.2. Estimation of the variance.

An estimation of the variance $\sigma$ is given by

$$\hat{\sigma}_i := \frac{\|R_{-i}^{-1/2} V_{-i}\|_2}{\sqrt{n-1}}.$$

Indeed, Eq. (9) gives that, under $\mathbb{H}_0$, it holds

$$\frac{\|R_{-i}^{-1/2} V_{-i}\|_2^2}{\sigma^2} \sim \chi^2(n-1),$$

where $\chi^2(n-1)$ is the chi-squared distribution with $n-1$ degree(s) of freedom. Since $V_{-i}$ is independent of $U_i$, note that $\varepsilon U_i$ and $\hat{\sigma}_i$ are independent. Furthermore, since $V_{-i}$ is Gaussian, remark that its norm and its direction are independent so that $V_{-i}/\hat{\sigma}_i$ and $\hat{\sigma}_i$ are independent. Recall that

$$\lambda_2^{i,\varepsilon} := \bigvee_{1 \leq j \neq i \leq p} \left\{ \frac{U_j - R_{ji} U_i}{1 - \varepsilon R_{ji}} \vee \frac{-U_j + R_{ji} U_i}{1 + \varepsilon R_{ji}} \right\},$$

and $V_{-i} := (U_1 - R_{1i}U_i, \ldots, U_{i-1} - R_{(i-1)i}U_i, U_{i+1} - R_{(i+1)i}U_i, \ldots, U_p - R_{pi}U_i)$. Eventually, remark that

$$(10) \qquad \varepsilon U_i, \ \frac{\lambda_2^{i,\varepsilon}}{\hat{\sigma}_i} \text{ and } \hat{\sigma}_i \text{ are mutually independent.}$$

## 3.3. Distribution of the test statistic.

Let $(i,\varepsilon)$ be in $[\![1,p]\!] \times \{\pm 1\}$. Recall that $\varepsilon U_i$ and $\lambda_2^{i,\varepsilon}$ are independent. In view of (10), observe that

$$T_1^{i,\varepsilon} := \frac{\varepsilon U_i}{\hat{\sigma}_i} \quad \text{and} \quad T_2^{i,\varepsilon} := \frac{\lambda_2^{i,\varepsilon}}{\hat{\sigma}_i},$$

are independent and, under $\mathbb{H}_0$, the random variable $T_1^{i,\varepsilon}$ a Student random variable with $n-1$ degree(s) of freedom. Define $(T_1, T_2)$ as

$$(11) \qquad (T_1, T_2) = \sum_{i=1}^{p} \sum_{\varepsilon=\pm 1} \left(\frac{\varepsilon U_i}{\hat{\sigma}_i}, \frac{\lambda_2^{i,\varepsilon}}{\hat{\sigma}_i}\right) \mathbb{1}_{\{\varepsilon U_i > \lambda_2^{i,\varepsilon}\}},$$

and recall that the events $\{\varepsilon U_i = \|U\|_\infty\} = \{\varepsilon U_i > \lambda_2^{i,\varepsilon}\}$ are almost surely disjoints.

**Lemma 3.** *Under $\mathbb{H}_0$, for each $i \in [\![1,p]\!]$ and $\varepsilon = \pm 1$, the random variable $T_2^{i,\varepsilon} = \lambda_2^{i,\varepsilon}/\hat{\sigma}_i$ has a density $p_{T_2^{i,\varepsilon}}^0$. Under $\mathbb{H}_0$, the joint density of $(T_1, T_2)$ is given by*

$$(12) \qquad \forall (\ell_1, \ell_2) \in \mathbb{R}^2, \quad p_{(T_1,T_2)}^0(t_1, t_2) = \mathbb{1}_{\{0 \leq t_2 \leq t_1\}} t_{n-1}(t_1) \sum_{i=1}^{p} \sum_{\varepsilon=\pm 1} p_{T_2^{i,\varepsilon}}^0(t_2),$$

where $t_{n-1}$ denotes the probability density function of the $t$-distribution with $n-1$ degree(s) of freedom.

*Proof.* One can check that $T_2^{i,\varepsilon}$ has a density, the reader may also consult Ylvisaker's theorem, see Theorem 1.22 in [AW09] for example. Observe that for all $(i,\varepsilon) \in [\![1,p]\!] \times \{\pm 1\}$, the random variable $T_2^{i,\varepsilon}$ is independent of the random variable $T_1^{i,\varepsilon}$. The result follows by (11). □

Recall that
$$T := \frac{1 - \mathbb{F}_{n-1}(T_1)}{1 - \mathbb{F}_{n-1}(T_2)},$$

where $\mathbb{F}_{n-1}$ denote the cumulative distribution function of the $t$-distribution with $n-1$ degrees of freedom. The expression of the joint density (12) shows that, conditionally to $T_2$, the random variable $T_1$ is distributed as a Student distribution conditioned to be greater than $T_2$. As a consequence, the conditional distribution of $T$ is uniformly distributed on $[0,1]$. We deduce that $T$ is uniformly distributed on $[0,1]$, as claimed.

## 4. THE TWO DIMENSIONAL CASE

In this section we assume that $p=2$ and $\Sigma = \mathrm{Id}_2$. Define

$$R = R(\rho) = \begin{pmatrix} 1 & \rho \\ \rho & 1 \end{pmatrix}, \tag{13}$$

with $\rho = \mathrm{Cov}(U_1, U_2) \in [-1, 1]$. Define the rejection region $\mathscr{R}_\alpha$ by

$$\mathbb{P}\{S \leq \alpha\} =: \mathbb{P}\{U = (U_1, U_2) \in \mathscr{R}_\alpha\}.$$

Note $\mathscr{R}_\alpha$ is symmetric about the origin and it is the non-convex disjoint union of four convex regions, namely

$$\mathscr{R}_\alpha^{+,1} = \{U_1 \geq \bar{\Phi}^{-1}(\alpha/2)\} \cap \{-g_\alpha(U_1)(1+\rho) + \rho U_1 \leq U_2 \leq g_\alpha(U_1)(1-\rho) + \rho U_1\},$$
$$\mathscr{R}_\alpha^{+,2} = \{U_2 \geq \bar{\Phi}^{-1}(\alpha/2)\} \cap \{-g_\alpha(U_2)(1+\rho) + \rho U_2 \leq U_1 \leq g_\alpha(U_2)(1-\rho) + \rho U_2\},$$
$$\mathscr{R}_\alpha^{-,1} = \{-U_1 \geq \bar{\Phi}^{-1}(\alpha/2)\} \cap \{-g_\alpha(-U_1)(1-\rho) + \rho U_1 \leq U_2 \leq g_\alpha(-U_1)(1+\rho) + \rho U_1\},$$
$$\mathscr{R}_\alpha^{-,2} = \{-U_2 \geq \bar{\Phi}^{-1}(\alpha/2)\} \cap \{-g_\alpha(-U_2)(1-\rho) + \rho U_2 \leq U_1 \leq g_\alpha(-U_2)(1+\rho) + \rho U_2\}.$$

where $g_\alpha(x) := \bar{\Phi}^{-1}(\bar{\Phi}(x)/\alpha) = \Phi^{-1}\left(1 - \frac{1-\Phi(x)}{\alpha}\right)$, for all $x \in \mathbb{R}$.

*Remark.* Observe this decomposition holds in any dimension $p$. The region $\mathscr{R}_\alpha$ given by $\mathbb{P}(S \leq \alpha) = \mathbb{P}(U \in \mathscr{R}_\alpha)$ is symmetric about the origin and is the non-convex disjoint union of $2^p$ convex regions of $\mathbb{R}^p$.

Note Anderson's inequality (see Lemma 4) is sufficient to establish the monotony of the Gaussian measure of symmetric convex set. Unfortunately, the region $\mathscr{R}_\alpha^c$ is not convex. However, when $p=2$, using an appropriate fibration, one may find a collection of sets (more general than $\mathscr{R}_\alpha^c$) satisfying a kind of generalization of Anderson's inequality. This is the object of the following proposition.

**Proposition 6.** *For $u \in \mathbb{R}$ and $\varepsilon = \pm 1$, define the non-centered diagonals:*

$$\Delta_u^\varepsilon = \{(x,y) \in \mathbb{R}^2 \; ; \; y = \varepsilon x + u\}.$$

*Let $\mathscr{T}$ be a set of $\mathbb{R}^2$ which is symmetric with respect to the two diagonals $\Delta_0^{+1}$ and $\Delta_0^{-1}$ and satisfies for all $u$ and $\varepsilon$:*

$$\mathscr{T} \cap \Delta_u^\varepsilon \text{ is an interval.} \tag{14}$$

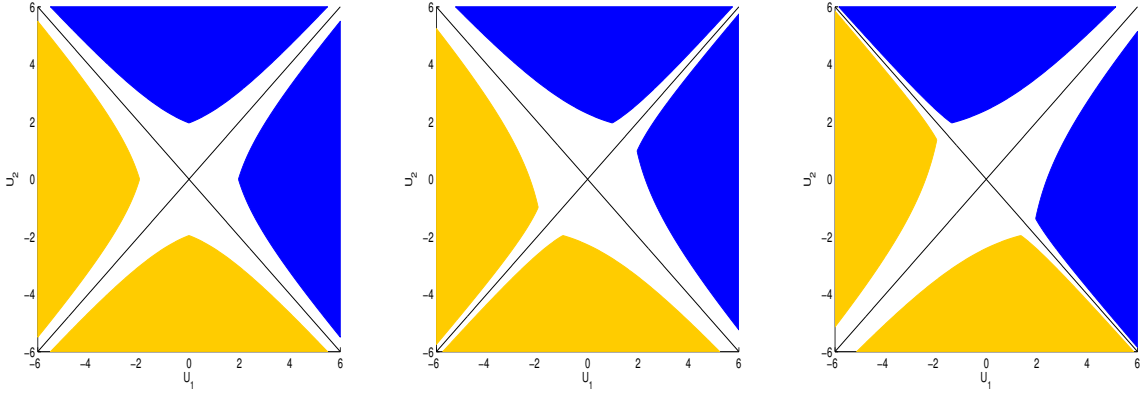

FIGURE 3. An illustration of region $\mathscr{R}_{0.05}^{+,1} \cup \mathscr{R}_{0.05}^{+,2}$ in (dark) blue and region $\mathscr{R}_{0.05}^{-,1} \cup \mathscr{R}_{0.05}^{-,2}$ in (bright) yellow for $\rho = 0$, $\rho = 0.5$ and $\rho = -0.7$ (from left to right).

*Set*

$$\Psi(\mu_1, \mu_2) = \mathbb{P}\{\mathcal{N}((\mu_1, \mu_2), R(\rho)) \in \mathcal{T}\}.$$

*Then the function $\Psi(u, v)$, which is obviously symmetric with respect to centered diagonals $\Delta_0^{+1}$ and $\Delta_0^{-1}$, is non-increasing along every half diagonal. More precisely, for every $u \in \mathbb{R}$ and $z \geq 0$, the two functions*

(15) $\qquad z \mapsto \Psi(u/2 - z, u/2 + z) = \Psi(u/2 + z, u/2 - z),$

(16) $\qquad z \mapsto \Psi(u/2 + z, -u/2 + z) = \Psi(u/2 - z, -u/2 - z),$

*are non-increasing.*

*Remark.* Note the result remains true if we multiply the matrix $R(\rho)$ by a scalar.

*Proof.* Set $R = R(\rho)$ for short. The proof relies on the fact that the eigenvectors of a two dimensional correlation matrix (such a $R$) are fixed and coincide with the diagonals of $\mathbb{R}^2$. Symmetry of $\mathcal{T}$ with respect to these diagonals is a key point in the proof.
First, we can uses a $\pi/4$ rotation and consider a variance-covariance matrix $\bar{R}$ which is diagonal and a set $\bar{\mathcal{T}}$ which is symmetric with respect to the two axes and whose intersections with a line parallel to one axe is an interval.
Let $\bar{\Psi}$ the expression of the function $\Psi$ after this rotation. We have to prove that $\bar{\Psi}$ is non-increasing along the relevant half-axes. Fix, for example, $\mu_1$ and consider, for $\mu_2 > 0$ the function

$$\mu_2 \mapsto \bar{\Psi}(\mu) = \mathbb{P}(\mathcal{N}(\mu, \bar{R}) \in \bar{\mathcal{T}}) = \int_{-\infty}^{+\infty} \varphi\left(\frac{u - \mu_1}{\sigma_1}\right) \mathbb{P}(\mathcal{N}(\mu_2, \sigma_2^2) \in I_u) du,$$

where $\sigma_1^2, \sigma_2^2$ are the diagonal elements of $\bar{R}$, $\mu = (\mu_1, \mu_2)$ and

$$I_u = \{v \in \mathbb{R} : (u, v) \in \bar{\mathcal{T}}\}.$$

Our hypotheses imply that for all $u$, $I_u$ is an interval that is symmetrical with respect to zero. Anderson's inequality (Lemma 4) implies directly that the Gaussian measure of $I_u$ is non-increasing as a function of $\mu_2$ so the function $\mu_2 \mapsto \bar{\Psi}(\mu_1, \mu_2)$ is non-increasing. This gives half of the statement, the other statement is obtained exactly in the same fashion by exchanging the roles of $\mu_1$ and $\mu_2$. $\qquad\square$

**Corollary 2.** $\forall \alpha \in (0,1)$, $\mathcal{R}_\alpha^c$ *satisfies the hypothesis of Proposition 6 so the power of the Spacing test for LARS is non-decreasing along the diagonals $\Delta_u^\varepsilon$ in the sense that is has exactly the same properties as those of the function $\Psi$ given by (15) and (16). In particular:*

- *Spacing test for LARS is unbiased,*
- *For each $\mu \in \mathbb{R}^2$, the function $t \mapsto \mathbb{P}_{t\mu}(S \leq \alpha)$ is non-increasing for $t \geq 0$.*

*Proof.* **Step 1**: If $(U_1, U_2)$ has distribution $\mathcal{N}(0, R)$ then it is is also the case of $(U_1, -U_2)$, $(-U_1, U_2)$ or $(-U_1, -U_2)$. This implies that $\mathcal{R}_\alpha^c$, which is computed under the null hypothesis, has the two required symmetry properties of Proposition 6.

**Step 2**: We consider now hypothesis (14). Consider $(u_1, u_2) \in \mathcal{R}_\alpha^{+,1}$. By definition of this region, it holds

$$u_1 \geq \bar{\Phi}^{-1}(\alpha/2),$$

and $\quad g_\alpha(u_1)(1-\rho) + \rho u_1 \geq u_2 \geq -g_\alpha(u_1)(1+\rho) + \rho u_1.$

Let $r > 0$ and consider the points $(u_1 + r, u_2 + r)$ and $(u_1 + r, u_2 - r)$. It is proven in Lemma 6 in the appendix that

$$g_\alpha(u_1 + r) \geq g_\alpha(u_1) + r.$$

As a consequence, for example, $g_\alpha(u_1+r)(1-\rho) + \rho(u_1+r) \geq g_\alpha(u_1)(1-\rho) + \rho(u_1) + r$ and this implies directly that $(u_1 + r, u_2 + r)$ and $(u_1 + r, u_2 - r)$ belong to $\mathcal{R}_\alpha^{+,1}$. The intersections of $\mathcal{R}_\alpha^{+,1}$ with the diagonals $\Delta_u^\varepsilon$ are half lines or empty sets. We have the same results for the three other regions in the same fashion and this implies that the intersections of $\mathcal{R}_\alpha^c$ with the diagonals $\Delta_u^\varepsilon$ are intervals. Finally, this result is true in $\beta$ because $R$ preserves symmetry properties along the diagonals. □

*Remark.* In dimension two, note

$$\Delta_0^1 \cup \Delta_0^{-1} = \{(U_1, U_2) \, ; \quad S(U_1, U_2) = 1\}.$$

In higher dimension, one has

$$\{(U_1, \ldots, U_p) \, ; \, S(U_1, \ldots, U_p) = 1\} = \bigcup_{i=1}^p \bigcup_{\varepsilon = \pm 1} \{\lambda_1 = \varepsilon U_i = \max_{j \neq i} |U_j|\}.$$

Observe the aforementioned set is not a hyperplane and so no orthogonal symmetry appears. The proof given in this section cannot be generalized to higher dimensions.

## 5. Numerical experiments on the power

**5.1. A Matlab Toolbox.** To compute the power of Spacing test for LARS using (5), we need to perform integration in high dimension. First, observe that (5) can be expressed as $n$-dimensional Gaussian integral (recall that $n$ is the rank of $X$). Indeed,

$$(17) \qquad \mathbb{P}_{\mu^\star}\{S \leq \alpha\} = \alpha \mathbb{E}_{\mu^\star}(W(U_1, \ldots, U_n))$$

where

$$(18) \qquad W(U_1, \ldots, U_n) = \sum_{i=1}^p \sum_{\varepsilon = \pm 1} \exp(\varepsilon \mu_i h_\alpha(\varepsilon U_i)) \mathbb{1}_{\mathcal{C}_{i,\varepsilon}}.$$

The aforementioned formula is a high dimensional Gaussian integral and we use a very efficient algorithm from A. Genz [Gen92, AG13], based on a reduction of the integral on the hypercube $[0,1]^n$ and Monte-Carlo Quasi Monte-Carlo (MCQMC) integration. In this fashion, Matlab programs *qsimvn* and *qsimvnef* provide powerful and robust numerical integration algorithms. The MCQMC routine is based on Kronecker or lattice sequences to compute integrals. In a second step, a Monte-Carlo (MC) layer is added to ensure unbiasedness and to compute the precision. Eventually, the QMC step is nested in the MC step in order to improve the speed of convergence, see [NC06] for example. A Matlab toolbox computing the power of Spacing test for LARS and based on Genz' routines is available on S. Mourareau's website [Mou15]. In addition, some practical examples are given.

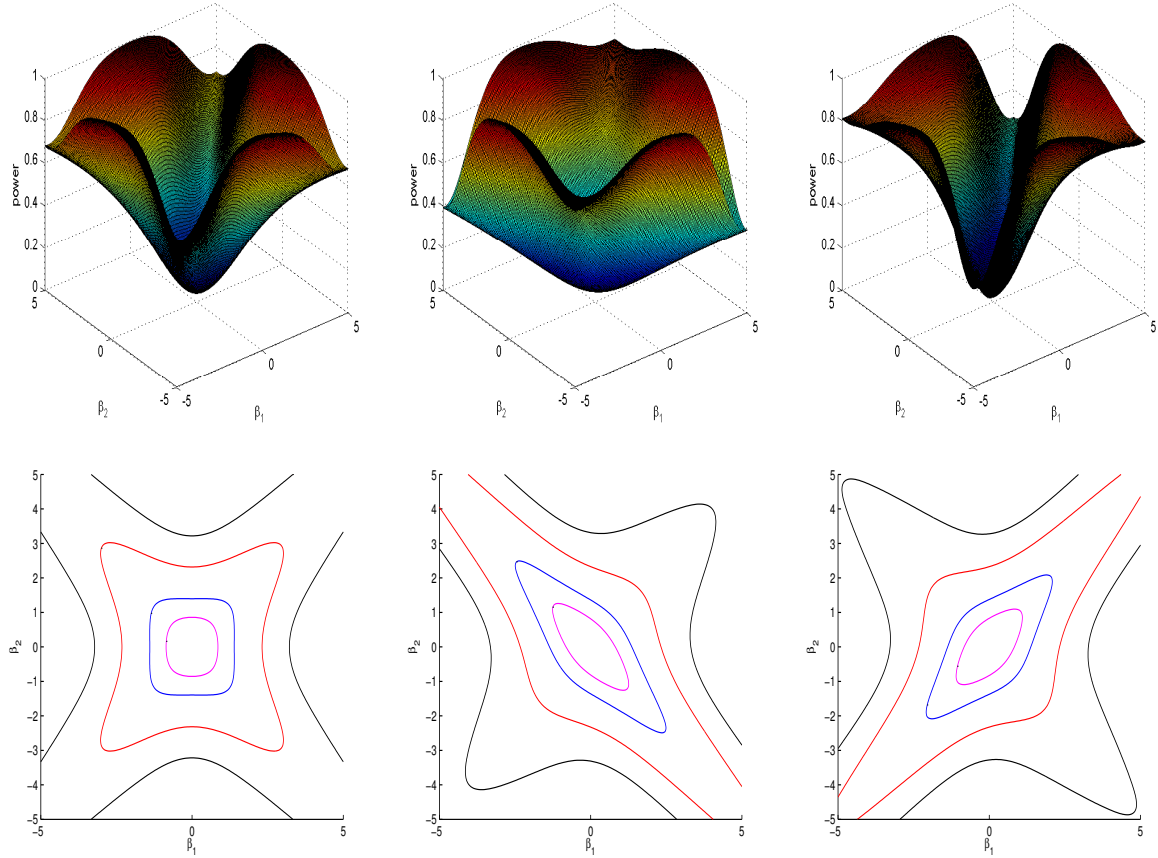

FIGURE 4. At the top, from left to right, power function $\beta \mapsto k_{\alpha,\rho}(\beta)$ for a significance level $\alpha = 0.05$ and correlations $\rho = 0$, $\rho = 0.5$ and $\rho = -0.4$. At the bottom, corresponding level sets of the power function, $k_{\alpha,\rho}(.) = 0.10, 0.20, 0.40$ and $0.70$.

5.2. **The two-dimensional case.** In dimension two, the power of Spacing test for LARS can be easily computed using numerical integration from (5). Consider the power function

$$k_{\alpha,\rho}(\beta) = \mathbb{P}(\mathcal{N}(R\beta, R) \in \mathcal{R}_\alpha),$$

where $R = R(\rho)$ is given by (13) and the region $\mathcal{R}_\alpha$ is defined in Section 4. The aforementioned power function is monotone in $\beta$ along the directions defined in Section 4 and it can be seen on Figure 4 that the variation of the power is minimal along the diagonal associated to the minimal eigenvalue of $R(\rho)$, see also Corollary 2.

5.3. **Pearson's chi-squared test versus the Spacing test for LARS.** We consider the standard goodness of fit test of the hypothesis

$$\mathbb{H}_0 : \text{``}\beta^\star \in \ker(X)\text{''} \quad \text{against} \quad \mathbb{H}_1 : \text{``}\beta^\star \notin \ker(X)\text{''}.$$

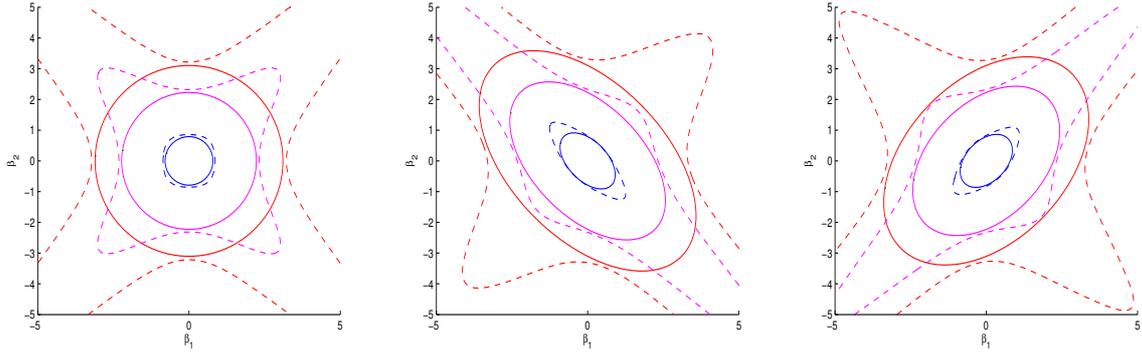

FIGURE 5. From left to right, level sets of the power functions of $S$ (dashed lines) and $T$ (plain line) for $\alpha = 0.05$ and $\rho = 0$, $\rho = 0.5$ and $\rho = -0.4$. We observe that the hypograph of the power function of $S$ is included in the corresponding one of $T$.

This test is defined by the statistic

$$T = \|Y\|_2^2$$

that follows a $\chi^2(n, \|X\beta\|_2^2)$ distribution where $\chi^2(a, b)$ denotes the $\chi^2$ distribution with $a$ degrees of freedom and non-centrality parameter $b$. Our aim is to compare this standard test with the Spacing test for LARS in different cases.

5.3.1. *The two-dimensional case.* In dimension two, considering the full model (s,n,p) = (2,2,2), we present a comparison of level sets of power functions for Spacing test for LARS and Pearson's chi-squared test, see Figure 5. It may suggest, from the comparison of level sets, that Pearson's chi-squared test is uniformly more powerful than Spacing test for LARS in the two-dimensional case.

5.3.2. *Higher dimensions.* In higher dimension, our experiments have the following frame. The design matrix $X$ is drawn from $n \times p$ independent standard Gaussian distribution. The target $\beta^\star$ has $s$ non zero entries independently and identically drawn from centered Gaussian distribution having variance 2 ("*large mean case*"), or variance 1 ("*medium mean case*") or from uniform random distribution on $[0,1]$ ("*small mean case*"). The choice of $s$, $n$ and $p$ concerns "*full*" models ($s = n = p$, see Figure 7), "*sparse*" models (see Figure 6) or "*very sparse*" models ($s \ll p$, see Figure 8).

5.3.3. *Conclusions.* Figure 5 suggests, from the comparison of level sets, that the $\chi^2$ test is uniformly more powerful than the Spacing test for LARS in the two dimensional case. Results from Figure 7 seem to confirm the interest of $\chi^2$ test in full models. However, Spacing test for LARS seems much more efficient in very high dimension cases when the signal presents a major gap between the dominant component and the rest (see Figure 8). When $\beta_i$ are of the same order of magnitude (meaning drawn from the same law), even in case of sparsity, $\chi^2$ test seems to be more powerful with respect to Spacing test for LARS (see Figure 6).

Simulations have been conducted with $\alpha = 0.05$ which is a classical choice. Other simulations with $\alpha = 0.01$ gave similar results. Finally, simulations involving higher level of $\alpha$ gave less marked results.

**Acknowledgment.** We thank Pr. Franck Barthe for valuable discussions.

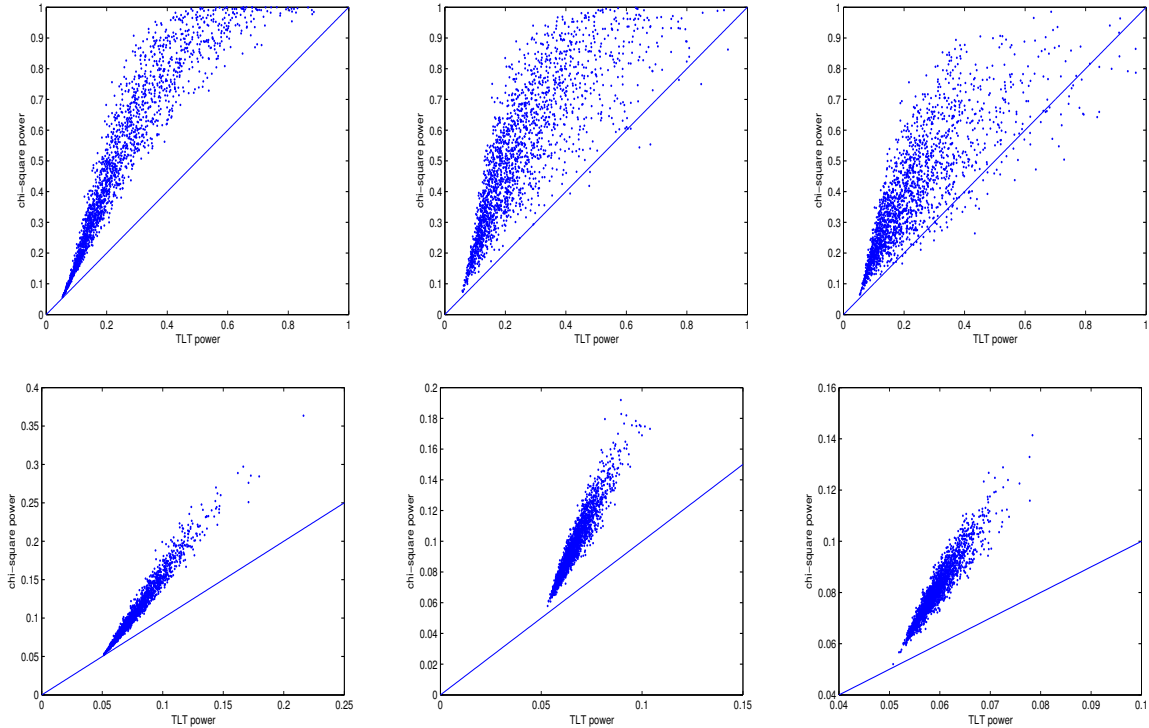

FIGURE 6. From left to right, 2.000 simulations of Spacing test for LARS's power versus $\chi^2$ power in various sparse cases $(s, n, p) = (5, 10, 50)$, $(10, 50, 100)$ and $(10, 100, 200)$. At the top, the mean $\beta$ is "large", while, at the bottom, the mean is "small" (see Section 5.3.2 for a definition). In both case, Pearson's chi-squared test seems more powerful in respectively 95, 94 and 99% of cases (large mean) and 91, 98 and 99% of cases.

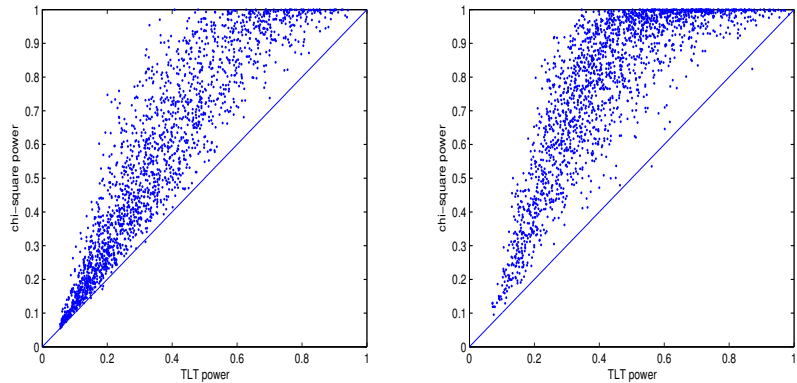

FIGURE 7. From left to right, 2000 simulations of Spacing test for LARS's power versus $\chi^2$ power in the "full" case $(s, n, p) = (5, 5, 5)$ and $(10, 10, 10)$ for a mixture of small, medium and high mean. As in dimension two, the $\chi^2$ test seems to give an improvement with respect to the Spacing test for LARS.

## APPENDIX

**Lemma 4** (Anderson's inequality for Gaussian measure [And55]). *Let $E$ be a convex set in $\mathbb{R}^p$, symmetric around the origin, and let $Z \sim \mathcal{N}_p(0, V)$. For all $t \geq 0$ and $\mu \in \mathbb{R}^p$*

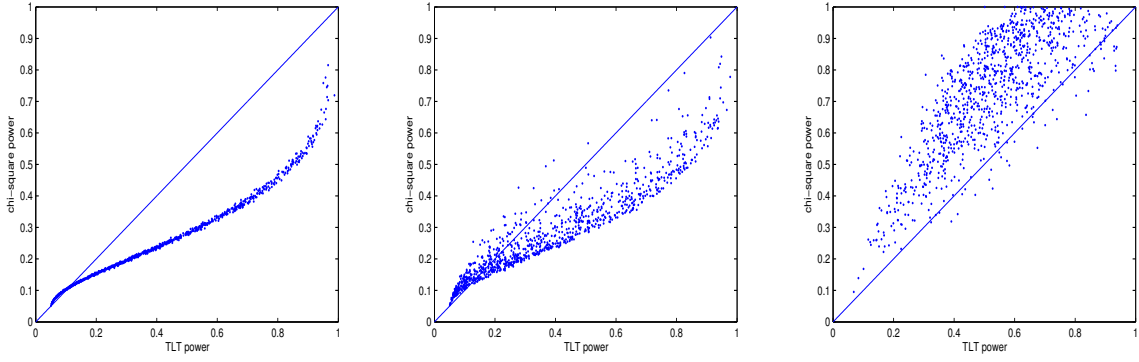

FIGURE 8. In the first instance (left), $(s,n,p) = (1,100,400)$ and the mean is drawn from $\mathcal{N}(\sqrt{2\log(p)},1)$. In the second one (center), one mean is drawn from $\mathcal{N}(\sqrt{2\log(p)},1)$ and others from $\mathcal{N}(0,1)$. In the third one (right), $(s,n,p) = (3,100,400)$ and all means are drawn from $\mathcal{N}(\sqrt{2\log(p)},1)$. When one mean is dominant, as in the first two cases, the Spacing test for LARS seems to be more efficient. However, when the difference between the two dominant means isn't large enough, the $\chi^2$ test seems to be more efficient.

*define*
$$\gamma_{E,\mu}(t) := \mathbb{P}(Z + t\mu \in E).$$
*Then $t \mapsto \gamma_{E,\mu}(t)$ is a non-increasing function.*

**Lemma 5.** *For all $\alpha \in [0,1]$ and for all $u \in \mathbb{R}$ such that $u \geq \bar{\Phi}(\alpha/2)$, it holds that the function $h_\alpha$ defined by (4) enjoys $h_\alpha \geq 0$, $h_\alpha$ is non-increasing, and $h_\alpha$ goes to zero at infinity.*

*Proof.* First, note $\alpha\bar{\Phi}(u) \leq \bar{\Phi}(u)$ and $\bar{\Phi}$ is non-increasing, to get that $h_\alpha \geq 0$. Compute the derivative and use that $h_\alpha \geq 0$ to show that $h'_\alpha(u) \leq \alpha - 1 \leq 0$. Eventually, we get that
$$\forall \alpha \in ]0,1], \ \forall u \geq u_0, \quad \bar{\Phi}(u+a) \leq \alpha\bar{\Phi}(u)$$
As $\bar{\Phi}$ is non-increasing, it implies that $u + a \geq \bar{\Phi}^{-1}(\alpha\bar{\Phi}(u))$ so that $a \geq h_\alpha(u) \geq 0$ which concludes the proof. □

**Lemma 6.** *$\forall \alpha \in [0,1], \forall u \geq \bar{\Phi}^{-1}(\alpha/2), \forall v \geq 0$, it holds*
$$(19) \qquad g_\alpha(u+v) \geq g_\alpha(u) + v,$$
*where $g_\alpha(u) = \Phi^{-1}\left(1 - \frac{1-\Phi(u)}{\alpha}\right) = \bar{\Phi}^{-1}(\bar{\Phi}(u)/\alpha)$.*

*Proof.* To prove (19), we show
$$\frac{\partial g_\alpha}{\partial u}(u,\alpha) = \frac{\varphi(u)}{\alpha\varphi(g_\alpha(u))} =: \frac{\varphi(u)}{j(u,\alpha)} \geq 1.$$
Use the fact that $\forall u \in \mathbb{R}, \ \varphi(u) \geq u\bar{\Phi}(u)$ to compute
$$\frac{\partial j}{\partial \alpha}(u,\alpha) = \varphi(g_\alpha(u)) - \frac{\bar{\Phi}(u)}{\alpha}g_\alpha(u) \geq 0$$
so
$$\frac{\partial g_\alpha}{\partial u}(u,\alpha) \geq \frac{\varphi(u)}{j(u,1)} = 1,$$
as claimed. □

# REFERENCES

[AG13]   J.-M. Azaïs and A. Genz, *Computation of the distribution of the maximum of stationary gaussian processes*, Methodology and Computing in Applied Probability **15** (2013), no. 4, 969–985.

[And55]  T. W. Anderson, *The integral of a symmetric unimodal function over a symmetric convex set and some probability inequalities*, Proceedings of the American Mathematical Society **6** (1955), no. 2, 170–176.

[AW09]   J.-M. Azaïs and M. Wschebor, *Level sets and extrema of random processes and fields*, John Wiley and Sons, 2009.

[BLPR11] K. Bertin, E. Le Pennec, and V. Rivoirard, *Adaptive dantzig density estimation*, Annales de l'IHP, Probabilités et Statistiques **47** (2011), no. 1, 43–74.

[BMvdG14] P. Bühlmann, L. Meier, and S. A. van de Geer, *Discussion: "a significance test for the lasso"*, Ann. Statist. **42** (2014), no. 2, 469–477.

[BRT09]  P. J. Bickel, Y. Ritov, and A. B. Tsybakov, *Simultaneous analysis of lasso and Dantzig selector*, Ann. Statist. **37** (2009), no. 4, 1705–1732. MR 2533469 (2010j:62118)

[BVDG11] P. Bühlmann and S. Van De Geer, *Statistics for high-dimensional data: Methods, theory and applications*, Springer, 2011.

[CDS98]  S. S. Chen, D. L. Donoho, and M. A. Saunders, *Atomic decomposition by basis pursuit*, SIAM J. Sci. Comput. **20** (1998), no. 1, 33–61. MR 1639094 (99h:94013)

[CT06]   E. J. Candès and T. Tao, *Near-optimal signal recovery from random projections: universal encoding strategies?*, IEEE Trans. Inform. Theory **52** (2006), no. 12, 5406–5425. MR 2300700 (2008c:94009)

[CT07]   ______, *The Dantzig selector: statistical estimation when p is much larger than n*, Ann. Statist. **35** (2007), no. 6, 2313–2351. MR 2382644 (2009b:62016)

[DC13]   Y. De Castro, *A remark on the lasso and the dantzig selector*, Statistics and Probability Letters **83** (2013), no. 1, 304 – 314.

[EHJT04] B. Efron, T. Hastie, I. Johnstone, and R. Tibshirani, *Least angle regression*, Ann. Statist. **32** (2004), no. 2, 407–499, With discussion, and a rejoinder by the authors. MR 2060166 (2005d:62116)

[Fuc05]  J.-J. Fuchs, *Recovery of exact sparse representations in the presence of bounded noise*, IEEE Transactions on Information Theory **51** (2005), no. 10, 3601–3608.

[Gen92]  A. Genz, *Numerical computation of multivariate normal probabilities*, Journal of computational and graphical statistics **1** (1992), no. 2, 141–149.

[HTF09]  T. Hastie, R. Tibshirani, and J. Friedman, *The elements of statistical learning*, vol. 2, Springer, 2009.

[HTW15]  T. Hastie, R. Tibshirani, and M. Wainwright, *Statistical learning with sparsity: The lasso and generalizations*, CRC Press, 2015.

[JN11]   A. Juditsky and A. Nemirovski, *Accuracy guarantees for l1-recovery*, Information Theory, IEEE Transactions on **57** (2011), no. 12, 7818–7839.

[LSST13] J. D. Lee, D. L. Sun, Y. Sun, and J. E. Taylor, *Exact post-selection inference with the lasso*, arXiv preprint arXiv:1311.6238 (2013).

[LTTT14a] R. Lockhart, J. Taylor, R. Tibshirani, and R. J. Tibshirani, *A significance test for the lasso*, The Annals of Statistics **42** (2014), no. 2, 413–468.

[LTTT14b] R. Lockhart, J. Taylor, R. J. Tibshirani, and R. Tibshirani, *Correction to rejoinder to "a significance test for the lasso"*, Ann. Statist. **42** (2014), no. 5, 2138–2139.

[LTTT14c] ______, *Rejoinder: "a significance test for the lasso"*, Ann. Statist. **42** (2014), no. 2, 518–531.

[Mou15]  S. Mourareau, http://www.math.univ-toulouse.fr/~smourare/power_test.m, March 2015.

[NC06]   D. Nuyens and R. Cools, *Fast algorithms for component-by-component construction of rank-1 lattice rules in shift-invariant reproducing kernel hilbert spaces*, Mathematics of Computation **75** (2006), no. 254, 903–920.

[Tib96]  R. Tibshirani, *Regression shrinkage and selection via the lasso*, J. Roy. Statist. Soc. Ser. B **58** (1996), no. 1, 267–288. MR MR1379242 (96j:62134)

[TLT14]  J. Taylor, J. Loftus, and R. J. Tibshirani, *Tests in adaptive regression via the kac-rice formula*, arXiv preprint arXiv:1308.3020v3 (2014).

[TLTT14] J. E. Taylor, R. Lockhart, R. J. Tibshirani, and R. Tibshirani, *Exact post-selection inference for forward stepwise and least angle regression*, arXiv preprint arXiv:1401.3889 (2014).

[vdGB09] S. A. van de Geer and P. Bühlmann, *On the conditions used to prove oracle results for the Lasso*, Electron. J. Stat. **3** (2009), 1360–1392. MR 2576316 (2011c:62231)

YDC is with Laboratoire de Mathématiques d'Orsay, Univ. Paris-Sud, CNRS, Université Paris-Saclay, 91405 Orsay, France.
  *E-mail address*: `yohann.decastro@math.u-psud.fr`

JMA and SM are with the Institut de Mathématiques de Toulouse (CNRS UMR 5219). Université Paul Sabatier, 118 route de Narbonne, 31062 Toulouse, France.
  *E-mail address*: `jean-marc.azais/stephane.mourareau@math.univ-toulouse.fr`